\documentclass{amsart}

\usepackage{enumerate}

\usepackage{enumitem}
\usepackage{graphicx}

\newtheorem{theorem}{Theorem}[section] 

\theoremstyle{definition}

\newtheorem{lemma}[theorem]{Lemma}
\newtheorem{conjecture}[theorem]{Conjecture}

\usepackage{hyperref}
\usepackage{breakurl}

\theoremstyle{remark}

\usepackage{algorithm, algpseudocode}

\begin{document}

\title[Computing isolated coefficients of the $j$-function]{Computing isolated coefficients of the $j$-function}


\author{Fredrik Johansson}
\address{Inria Bordeaux-Sud-Ouest and Institut Math.\ Bordeaux,  U. Bordeaux, \newline33400 Talence, France}
\curraddr{}
\email{fredrik.johansson@gmail.com}
\thanks{}



\date{}

\dedicatory{}

\begin{abstract}
We consider the problem of efficiently computing isolated coefficients $c_n$ in the
Fourier series of the elliptic modular function $j(\tau)$. We show that
a hybrid numerical-modular method with complexity $n^{1+o(1)}$
is efficient in practice.
As an application, we locate
the first few values of $c_n$ that are prime, the first occurring at $n = 457871$.
\end{abstract}

\maketitle

\section{Introduction}

The coefficients $c_n$ appearing in the Fourier series
of the elliptic modular function
$j(\tau) = \sum_{n=-1}^{\infty} c_n e^{2\pi i n \tau}$
have many remarkable arithmetical properties.
This integer sequence (A000521 in Sloane's OEIS~\cite{OEIS000521}) begins
$$c_{-1} = 1, \quad c_0 = 744, \quad c_1 = 196884, \quad c_2 = 21493760, \quad \ldots$$
and is perhaps most notorious
for its connection with the monster group (a correspondence known as \emph{monstrous moonshine} \cite{conway1979monstrous}).

The numerical growth rate of these coefficients is known:
$c_n$ has around $\beta \sqrt{n}$ digits where $\beta = 4 \pi / \log(10) \approx 5.46$.
More precisely, explicit lower and upper bounds for $n \ge 1$ are given by \cite[Theorem~1.1]{Brisebarre2005}
$$c_n = \frac{e^{4 \pi \sqrt{n}}}{\sqrt{2} n^{3/4}} \left(1 - \frac{3}{32 \pi \sqrt{n}} + \varepsilon_n \right), \quad |\varepsilon_n| \le \frac{0.055}{n}.$$

Our concern in this work will be the exact calculation
of $c_n$ for large $n$.
Baier and K\"{o}hler~\cite{Baier2003} survey several strategies,
concluding that a formula by Zagier and Kaneko
is the ``most efficient method''.
However, their analysis only considers
recursive calculation which does not yield the best
possible complexity. Using fast power series
arithmetic, applicable to many of the formulas
discussed by Baier and K\"{o}hler, it is possible to compute
$c_{-1},\ldots,c_n$ simultaneously in $n^{1.5+o(1)}$ time (bit operations),
or in $n^{1+o(1)}$ time modulo a fixed number~$M$.
These complexities are essentially optimal
since they are quasilinear in the size of the output.

A more challenging problem is to compute an isolated
value $c_n$ quickly.
It is presently unknown whether an algorithm with
quasioptimal
$n^{0.5+o(1)}$ time complexity exists,
but it is possible to achieve 
$n^{1+o(1)}$ time and $n^{0.5+o(1)}$ space complexity
using numerical evaluation of an infinite series for $c_n$ derived
by Petersson~\cite{petersson1932} and Rademacher~\cite{Rademacher1938},
improving on the
$n^{1.5+o(1)}$ time and $n^{1+o(1)}$ space
needed with the power series method.
Unfortunately, Rademacher, Baier and K\"{o}hler
and several authors have dismissed this method
as impractical for computations.
We analyze this method below and find that it is useful,
but that current error bounds are too pessimistic to make it practical
if we insist on rigorous results.

As a solution to this problem, we propose a rigorous hybrid algorithm
that uses the Petersson-Rademacher series
to obtain the high bits of $c_n$ together with a
power series calculation or a suitable congruence
to obtain the low bits.
This method requires $n^{1+o(1)}$ time and space,
or $n^{0.5+o(1)}$ space for $n$ of special form,
and is highly efficient in practice.

This work was prompted by a MathOverflow post by David Feldman
which asked whether the sequence $c_n$ contains prime numbers~\cite{MO2}.
We can answer this in the affirmative: an
exhaustive search of $n \le 2 \cdot 10^7$
finds seven prime values of $c_n$ in this range (Table~\ref{tab:primes}), the first
occurring at $n = 457871$.

\section{The power series method}

A natural way to compute $c_n$ is to
expand $j(\tau) = \sum_{k=-1}^{\infty} c_k q^k$ as a truncated formal power series
in $q$ and read off the last coefficient.
This method obviously gives all the coefficients up to $c_n$ if
we so prefer.

There are many ways to express $j(\tau)$ in terms of
simple functions suitable for power series computations.
The following formula is economical and works in any characteristic:
define the Eisenstein series
$E_4 = 1 + 240 \sum_{n=1}^{\infty} \sigma_3(n) q^n$
where $\sigma_x(n)=\sum_{d\mid n} d^x$ is the divisor function,
and denote by
$\phi = \sum_{n \in \mathbb{Z}} (-1)^n q^{n(3n-1)/2}$
the Dedekind eta function without
leading factor $q^{1/24}$.
Then
\begin{equation}
j(\tau) = \frac{1}{q} \left(\frac{E_4}{\phi^8}\right)^3.
\label{eq:je4formula}
\end{equation}

Noting that $\phi^4 = \phi^3 \cdot \phi$ can be generated
cheaply as a sparse product, the cost of evaluating \eqref{eq:je4formula} is essentially 2 multiplications, 2 squarings and 1 inversion of power series,
plus the construction of $E_4$ which can be done
in similar time to a multiplication using the sieve of Eratosthenes
(alternatively, $E_4$ can be constructed from Jacobi theta functions using further multiplications).
For a small speedup, we note that the final multiplication $(E_4 / \phi^8)^3 = (E_4 / \phi^8)^2 \cdot (E_4 / \phi^8)$
can be replaced by an linear summation if we only want the last coefficient.

The computations can be done in $n^{1.5+o(1)}$ time using FFT multiplication together with Newton inversion,
since we have power series of length $n$ with coefficients $n^{0.5+o(1)}$ bits in size.
The memory requirement is $n^{1.5+o(1)}$ bits if we work
directly over $\mathbb{Z}$. 
Computing modulo a $b$-bit integer $M$, the complexity for determining $c_n \bmod M$ reduces to $(nb)^{1+o(1)}$.
The memory requirement
for computing the single integer $c_n$
can thus be reduced to $n^{1+o(1)}$
if we compute separately modulo small pairwise coprime integers $m_k$ with $M = m_1, \ldots, m_k$, $c_n < M$,
and then reconstruct $c_n$ using the Chinese remainder theorem.

\section{The Petersson-Rademacher series}

Given any $n \ge 1$, the Petersson-Rademacher series for $c_n$ is the convergent series (letting $N \to \infty$)
\begin{equation}
c_n = \frac{2 \pi }{\sqrt{n}} \sum_{k=1}^N \frac{S(n,-1,k)}{k} \,I_1\!\left(\frac{4 \pi \sqrt{n}}{k}\right) + R_N(n)
\label{eq:radseries}
\end{equation}
in which $S(a,b;k)$ denotes a Kloosterman sum, $I_{\nu}(x)$ denotes a modified Bessel function of the first kind,
and the remainder term can be shown to satisfy
\begin{equation}
|R_N(n)| \le \frac{72 \pi}{\sqrt{n}} N^{3/4} \,I_1\!\left(\frac{4 \pi \sqrt{n}}{N}\right).
\label{eq:rbound}
\end{equation}

For a proof of the truncation bound as well as a generalization
to the corresponding coefficients for the function $j(\tau)^m$,
see Brisebarre and Philibert~\cite{Brisebarre2005}.

The Kloosterman sum is the exponential sum
$$S(a,b;k) = \sum_{\gcd(x,k) = 1} e^{2\pi i (ax+by) / k}$$
where the index $x$ ranges over $0 \le x < k$ and $y$ is any solution of $xy \equiv 1 \bmod k$.

The Petersson-Rademacher series for $c_n$ is analogous to the Hardy-Ramanujan-Rademacher
formula for the integer partition function $p(n)$,
which allows computing $p(n)$ in
quasioptimal time $n^{0.5+o(1)}$~\cite{Johansson2012hrr,Johansson2014thesis}.
The idea is that adding $\Theta(n^{0.5})$
terms of the series gives an approximation $y$ with $|y - p(n)| < 0.5$,
yielding the correct result when rounded to the nearest integer.
Although there are $\Theta(n^{0.5})$ terms and the result has
$\Theta(n^{0.5})$ bits, the overall complexity is $n^{0.5+o(1)}$ rather than $n^{1+o(1)}$
when the algorithm is implemented carefully
since the bit sizes of the terms fall off as a hyperbola after the first term.

By a similar analysis of the Petersson-Rademacher series for $c_n$, one can show the following:

\begin{theorem}
The integer $c_n$ can be computed using $n^{1+o(1)}$ bit operations and $n^{0.5+o(1)}$ bits of space.
\label{thm:complexity}
\end{theorem}

The reason why we do not get an $n^{0.5+o(1)}$ algorithm
is that computing the Kloosterman sum $S(a,b;k)$ ostensibly requires adding $O(k)$ terms.
The computation of 
$S(a,b;k)$ can be reduced to the computation of the shorter sums $S(a_i,b_i;q_i)$
for each prime power $q_i = p_i^{e_i}$ in the factorization of $k$,
and closed formulas are known when $e_i \ge 2$~\cite{whiteman1945note}, but for prime modulus
$q_i = p$ no algorithm better than $O(p)$ summation is currently known
(the existence of such an algorithm would immediately lead
to a better complexity bound for computing $c_n$).
The corresponding exponential sums in the series for $p(n)$
admit a complete factorization into simple trigonometric expressions
and can therefore be computed rapidly.

\subsection{Analysis of the error bound}

Although the Petersson-Rademacher series regrettably does not yield
an $n^{0.5+o(1)}$ algorithm with current technology,
the problem with the method is not the asymptotic
$n^{1+o(1)}$ complexity (which is quite serviceable),
but the hidden constant factors.

To make calculations explicit, we may combine \eqref{eq:rbound}
with the following bounds for the Bessel function $I_1(x)$, accurate when $x \to \infty$ and $x \to 0$ respectively
(obtained from the asymptotic expansion and the Taylor series at $x = 0$).

\begin{lemma}
For $x > 0$, $I_1(x) < e^x / \sqrt{2 \pi x}$.
\label{lem:besselexp}
\end{lemma}

\begin{lemma}
For $0 < x < 0.1$, $I_1(x) < 0.501x$.
\label{lem:bessel0}
\end{lemma}

A direct calculation using Lemma~\ref{lem:bessel0} gives, for instance:

\begin{theorem}
If $N \ge \max(C_0, C_1 \sqrt{n})$ where $C_0 = 6.7 \cdot 10^{13}$ and $C_1 = 40 \pi$, then $|R_N(n)| \le 0.499$.
\end{theorem}

This shows that $\Theta(n^{0.5})$ terms are sufficient to determine $c_n$ (as needed
in the proof of Theorem~\ref{thm:complexity}), noting that
bounding the truncation error by $0.499$ gives some wiggle room
for floating-point error in the approximation of the sum.

The constant $C_0$ makes the method virtually useless, as we would
have to perform some $C_0^2 \approx 10^{26}$ operations
to compute any $c_n$.
The constants $C_0$ and $C_1$ in this theorem are not optimal,
but $C_0$ cannot be brought below $16 \cdot (12 \pi)^8 \approx 6.53 \cdot 10^{13}$
using the bound \eqref{eq:rbound} (without a corresponding pessimistic increase of $C_1$).
Rademacher, working with a somewhat worse error bound than
\eqref{eq:rbound}, similarly concluded:

\begin{quotation}
``Unfortunately, the convergence of [...] is rather slow, so that we should need quite a number of
terms in order to get an error which is safely less than $1/2$.''
\end{quotation}

\begin{figure}[h]
\centering
\includegraphics[width=0.85\textwidth]{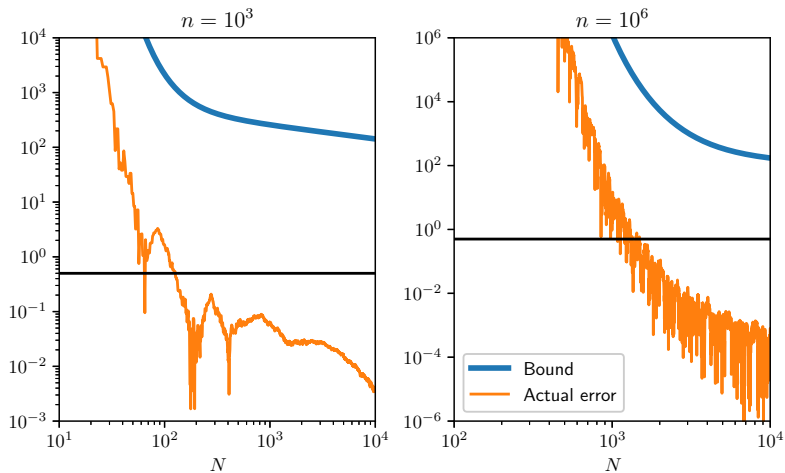}
\caption{The actual error $R_N(n)$ in the approximation of $c_n$ after summing $N$ terms of the Petersson-Rademacher series \eqref{eq:radseries}, compared with the bound~\eqref{eq:rbound}. The horizontal line locates $1/2$.}
\label{fig:errplot}
\end{figure}

How slow is the convergence really? We can run some experiments to get an idea.
Figure~\ref{fig:errplot} compares the actual rate of convergence of
the Petersson-Rademacher series with the bound \eqref{eq:rbound}.
It turns out that the bound is quite pessimistic,
and if we cut off the summation heuristically,
the computation becomes practical.
However, simply stopping when the partial sum seems to be
very close to an integer is dangerous. The erratic nature of the terms
(clearly visible in the figure, and observed in similar
contexts by other authors, e.g.~\cite{McLaughlin2012})
means that the sum can stabilize within $\varepsilon$ of an integer
for rather small $\varepsilon$ and remain there through many consecutive terms before suddenly
making a large jump. We could probably find a good empirical fit
for the true error, but we prefer a rigorous analysis.

This brings us to the question of how \eqref{eq:rbound} can be
improved.
The ingredients in the derivation of this bound are \cite[Section 5.1]{Brisebarre2005}:
\begin{itemize}
\item The Weil bound $|S(a,b;k)| \le \sqrt{\sigma_0(k)} \sqrt{\gcd(a,b,k)} \sqrt{k}$.
\item The bound $\sigma_0(k) \le 9k^{1/4}$ for the number of divisors of $k$.
\item Bounding the sum of terms in \eqref{eq:radseries} by a sum of absolute values of the terms.
\end{itemize}

The bound $\sigma_0(k) \le 9k^{1/4}$ is not optimal.
However, this bound is quite reasonable over the relevant non-asymptotic range of $k$.
If we could estimate $\sigma_0(k)$ by its \emph{average}
value $\log(k)$, and show that the deviations
from the average give a negligible contribution,
we would get a constant $C_0$ around
$3 \cdot 10^7$, which is much better but still rather impractical
for computations.

It therefore seems that a useful error bound
will require a much more involved analysis of cancellation in
sums of Kloosterman sums.
This is a well-studied problem~\cite{sarnak2009linnik},
but it appears that all published results stronger than
the Weil bound
are asymptotic without explicit constants.
We leave it as an open problem to prove sharp bounds for $R_N(n)$.

\section{Hybrid algorithm}

To get around the issues discussed in the previous
section without sacrificing rigor, we make the following
modification: instead of computing $c_n$ directly,
we assume that $c_n$ is known modulo some integer $M$.
If we denote by $r$ the unique residue of $c_n$ with $0 \le r < M$,
we have $c_n = M [(c_n - r)/M] + r$.
We can now compute $(c_n - r) / M$ using \eqref{eq:radseries},
stopping the summation when the bound for $|R_N(n)|$ is just less than
$M/2$ instead of $1/2$.
If we choose $M$ large enough,
the truncation error is in the exponentially-large domain of the Bessel function (Lemma~\ref{lem:besselexp})
where the available error bound works well.

For a general index $n$, we can compute the residue of $c_n$ in $\mathbb{Z}/M\mathbb{Z}$
using the power series \eqref{eq:je4formula}.
We can either compute a single power series over $\mathbb{Z}/M\mathbb{Z}$
or choose a composite modulus, say $M = p_1 p_2 \cdots p_k$,
and compute the value for each $\mathbb{Z}/p_i\mathbb{Z}$
in order to minimize memory usage.
If $M$ is bounded, the memory usage will be $n^{1+o(1)}$,
and the running time will be $n^{1+o(1)}$.
This is not as good as the Petersson-Rademacher series would be if we had
an optimal error bound, but it is better than computing $c_n$ using
power series alone.

The optimal choice of $M$ will depend on the
implementation-dependent
relative speeds of the power series arithmetic
and the numerical calculations for the Petersson-Rademacher series.
In our experiments, we have obtained the best performance
with $M$ slightly larger than $2^{10}$ for computing isolated values of $c_n$.
We get relatively uniform performance with $M$ anywhere between $2^{10}$
and $2^{20}$ while the running time increases sharply if $M$
is smaller than~$2^{9}$.
If we want to compute a range of values of $c_n$, it is more
efficient to choose a larger $M$, say word-size $M \approx 2^{64}$ (or bigger,
if sufficient memory is available),
precomputing the residues $c_{-1},\ldots,c_n$ once
before using the Petersson-Rademacher series to determine the high bits for each coefficient of interest.

\subsection{Special indices}

When $n$ has special form, we can use known congruences for $c_n$
to determine a residue. For $a > 1$, we have
\begin{equation*}
\begin{matrix}
c_{2^a m} & \equiv & -2^{3a+8} 3^{a-1} \sigma_7(m) & \pmod {2^{3a+13}}, & \quad m \text{ odd} \\
c_{3^a m} & \equiv & \mp 3^{2a+3} 10^{a-1} \sigma(m)/m & \pmod {3^{2a+6}}, & \quad m \equiv \pm 1 \pmod 3 \\
c_{5^a m} & \equiv & -5^{a+1} 3^{a-1} m \sigma(m) & \pmod {5^{a+2}} \\
c_{7^a m} & \equiv & -7^a 5^{a-1} m \sigma_3(m) & \pmod {7^{a+1}} \\
\end{matrix}
\end{equation*}
and there are also congruences modulo 11 and 13~\cite{aas1964congruences,Newman1958}.
If $n$ is divisible by sufficiently many small primes,
we can thus immediately construct a residue modulo an~$M$
that is large enough (say $M > 1000$) to use the Petersson-Rademacher series.
In that case, we do not need a power series evaluation,
and the entire computation can be done using $n^{0.5+o(1)}$ memory.
Indeed, it is sufficient that $n$ is even, since this gives the comfortably
large $M \ge 65536$.
If the $M$ obtained from congruences alone is too small, we can amend
it with a power series computation modulo some larger prime.

\section{Computations}

We have implemented evaluation of $c_n$ using direct power series
methods as well as the hybrid algorithm, using Flint~\cite{Har2010}
for power series and integer operations
and Arb~\cite{Joh2017} for arbitrary-precision ball arithmetic.
The following computations were performed on a laptop
with 16~GB RAM and an Intel i5-4300U CPU.

\begin{table}
\begin{center}
\caption{Time in seconds to compute $c_n$ or $c_{n'}$ where $n'$ = next prime after $n$. PS\ denotes
the power series method using arithmetic directly in $\mathbb{Z}$. PS2\ denotes the power series method
using a multimodular approach (conserving memory). NEW denotes the
hybrid algorithm using a congruence together with the Petersson-Rademacher series.
The entry (mem) indicates that the computation could not be run on the
machine due to insufficient memory.}
\label{tab:timings}
\renewcommand{\baselinestretch}{1.25}
\begin{small}
\begin{tabular}{ c c | c c c | c c } 
$n$  & Digits in $c_n$ & PS\ ($c_n$) & PS2 ($c_n$) & NEW ($c_n$) & $n'$ & NEW ($c_{n'}$)\\ \hline
$10^2$ & 53     & 0.00011  & 0.00017    &   0.00033 & $10^2 + 1$ & 0.00087 \\
$10^3$ & 171    & 0.0068   & 0.015      &   0.0016  & $10^3 + 9$ & 0.0041 \\
$10^4$ & 543    & 0.33     & 0.96       &   0.0054  & $10^4 + 7$ & 0.027 \\
$10^5$ & 1722   & 15       & 48         &   0.025   & $10^5 + 3$ & 0.22 \\
$10^6$ & 5453   & 625      & 2169       &   0.13    & $10^6 + 3$ & 2.5 \\
$10^7$ & 17253  & (mem)       & $\approx$ $8.2 \cdot 10^4$      &   0.83    & $10^7 + 19\!\!\!\!$ & 32 \\
$10^8$ & 54569  &          & $\approx$ $3.2\cdot 10^6$           &   6.5     & $10^8 + 7$ & 380 \\
$10^9$ & 172575 &          & (mem)        &   60      & $10^9 + 7$ & (mem) \\
$10^{10}$ & 545743  &      &            &   636     & & \\
\end{tabular}
\end{small}
\end{center}
\end{table}

\subsection{Large values}

As a benchmark (Table~\ref{tab:timings}), we compute $c_n$ for various powers of ten $n = 10^k$
as well as $c_{n'}$ where $n'$ is the first prime number
after $10^k$. Powers of ten are numbers of special form
($n$ being divisible by $2^k 5^k$),
meaning that the hybrid algorithm can use a congruence,
skipping the power series evaluation.
The subsequent prime number $n'$ is of generic form (worst-case input) for the algorithm,
forcing a power series evaluation to determine a residue.

For comparison purposes, we have implemented two versions of
the power series algorithm to compute $c_n$: the first
using arithmetic in $\mathbb{Z}$ and the second
using arithmetic modulo many small primes to conserve memory
(which turns out to be roughly three times slower).
The first power series implementation
is roughly equivalent to the function \texttt{j\_invariant\_qexp}
in SageMath~\cite{Sag2020}.

We see that the hybrid algorithm is faster
than the power series method already around $n = 10^3$.
At $n = 10^6$, it is 250 times faster for the generic input
(prime~$n$) and 4800 times faster
for the special-form input (power-of-ten $n$).
The $n^{1+o(1)}$ complexity of the hybrid
algorithm is apparent in the running times.\footnote{Arb does not presently
use a quasi-optimal algorithm for the Bessel function $I_1(x)$ at high precision,
so if we were to continue the table beyond $10^{10}$,
the timings would likely get worse.}

\subsection{Prime values}

The plethora of congruences satisfied by $c_n$
conspire to rule out $c_n$ being a prime number for small $n$,
but nothing suggests that this pattern must hold
asymptotically.
Indeed, while $c_0,\ldots,c_{70}$ are all
divisible by 2, 3 or 5, the smallest factor of $c_{71}$ is 353.

To search for prime values of $c_n$, we used a single power series
evaluation
to compute several $c_n$ simultaneously modulo the primorial
$M = 2 \cdot 3 \cdot \cdots 47$, which fits in a 64-bit word.
We then selected the entries with $\gcd(c_n, M) = 1$
and computed their full values 
using the hybrid algorithm with $c_n \bmod M$ as precomputed input.
Prime values of $c_n$ were then identified using a standard probabilistic
primality test (trial factoring to rule out simple composites
followed by the BPSW test).

To search up to $n = 2 \cdot 10^7$, the initial power series
calculation took three minutes. Filtering by $\gcd(c_n, M) = 1$
left 28971 candidate $c_n$ to compute exactly (around 7 hours) and
check for primality (around 45 hours).

\begin{table}
\begin{center}
\caption{All prime values of $c_n$ with $n \le 2 \cdot 10^7$.}
\label{tab:primes}
\renewcommand{\baselinestretch}{1.25}
\begin{small}
\begin{tabular}{ c | c | c } 
$n$  & Digits in $c_n$ & $c_n$ \\ \hline
457871 & 3689 & $3080163651 \ldots 2714076699$ \\
685031 & 4513 & $2912989222 \ldots 8765523019$ \\
1029071 & 5532 & $4025516131 \ldots 1099172019$ \\
1101431 & 5723 & $8315472348 \ldots 7940410921$ \\
9407831 & 16734 & $9603424490 \ldots 8550890201$ \\
11769911 & 18718 & $5971402918 \ldots 6331345197$ \\
18437999 & 23429 & $4474491259 \ldots 2242965811$ \\
\end{tabular}
\end{small}
\end{center}
\end{table}

Table~\ref{tab:primes} lists the first prime values of $c_n$.
At the time of writing, Jeremy Rouse reports having
certified the primality of $c_{457871}$ using ECPP.
The remaining numbers are only confirmed as probable primes,
but the BPSW test has no known counterexamples,
and we anticipate that ECPP certificates can
be generated with some months of computation.

Although primes seem to appear sparsely in the sequence $c_n$ (much more sparsely than for the
partition function $p(n)$, for example), we speculate:

\begin{conjecture}
There are infinitely many $n$ such that $c_n$ is prime.
\end{conjecture}

With more memory and a large number of cores, the search could easily be extended
further (at least to $n = 10^9$). The method presented here could also be used
to investigate other divisibility properties of the numbers $c_n$.

\section{Generalization}

The methods discussed here are not specific to the $j$-function.
The Petersson-Rademacher series can be generalized to
any modular function of weight 0 (see \cite[Theorem~5.1]{Brisebarre2005}),
and similar series can be constructed
for the coefficients of a wide range of modular forms,
or viewed combinatorially,
for various partition-type sequences~\cite{sills2010towards}.
An interesting problem is to automate the efficient
computation of such coefficients.
Obtaining tight truncation bounds for Rademacher-type series
seems difficult at the moment,
but with the hybrid
numerical-modular approach, more crude and generally
applicable bounds can be used.

\section{Source code}

The author has made all source code and data behind this
paper publicly available at \url{https://github.com/fredrik-johansson/jfunction}

\bibliographystyle{plain}
\bibliography{references}

\end{document}